\documentclass[oneside]{amsart}
\usepackage{longtable}
\usepackage{amssymb}

\makeatletter

\providecommand{\tabularnewline}{\\}

\theoremstyle{plain}
	\newtheorem{thm}{Theorem}[section]
\theoremstyle{plain}
	\newtheorem*{thm*}{Theorem}
\theoremstyle{plain}
	\newtheorem{lem}[thm]{Lemma}
\theoremstyle{remark}
	\newtheorem{rem}[thm]{Remark}
\theoremstyle{plain}
	\newtheorem{cor}[thm]{Corollary}
\theoremstyle{plain}
	\newtheorem{prop}[thm]{Proposition}
\theoremstyle{remark}
	
\theoremstyle{remark}
	\newtheorem*{acknowledgement*}{Acknowledgement}

\newcommand{\sgn}{\mathrm{sgn}}
\newcommand{\gl}{\mathrm{GL}}
\newcommand{\Ind}{\mathrm{Ind}}
\newcommand{\norm}{\mathrm{norm}}
\newcommand{\cfq}{\mathbb{F}_q}
\newcommand{\cfp}{\mathbb{F}_p}
\newcommand{\mdpd}[4]{\left(\begin{array}{cc}#1 & #2 \\ #3 & #4 \end{array}\right)}
\newcommand{\ydg}[1]{ {}_{#1}^{#1}\mathcal{YD}}
\newcommand{\mbq}{\mathbf{q}}
\newcommand{\nic}{\mathfrak{B}}
\newcommand{\glfq}{\mathbf{GL}(2,\mathbb{F}_q)}
\newcommand{\slfq}{\mathbf{SL}(2,\mathbb{F}_q)}

\newcommand{\pslfq}{\mathbf{PSL}(2,\mathbb{F}_q)}
\newcommand{\ZZ}{\mathbb{Z}}
\newcommand{\NN}{\mathbb{N}}
\newcommand{\ra}[1]{\mathcal{R}_{#1}}
\newcommand{\Dospuntos}[2]{%
	\setlength{\unitlength}{.12cm}
	\rule[-3\unitlength]{0pt}{8\unitlength}
	\begin{picture}(14,5)(0,3)
		\put(1,2){\circle{2}}
		\put(13,2){\circle{2}}
		\put(1,5){\makebox[0pt]{\scriptsize #1}}
		\put(13,5){\makebox[0pt]{\scriptsize #2}}
	\end{picture}
}
\newcommand{\Dchaintwol}[3]{%
	\setlength{\unitlength}{.12cm}
	\rule[-3\unitlength]{0pt}{8\unitlength}
	\begin{picture}(14,5)(0,3)
		\put(1,2){\circle{2}}
		\put(2,2){\line(1,0){10}}
		\put(13,2){\circle{2}}
		\put(1,5){\makebox[0pt]{\scriptsize #1}}
		\put(7,4){\makebox[0pt]{\scriptsize #2}}
		\put(13,5){\makebox[0pt]{\scriptsize #3}}
	\end{picture}
}
\newcommand{\Dchaintwo}[3]{%
	\setlength{\unitlength}{.08cm}
	\rule[-3\unitlength]{0pt}{8\unitlength}
	\begin{picture}(14,5)(0,3)
		\put(1,2){\circle{2}}
		\put(2,2){\line(1,0){10}}
		\put(13,2){\circle{2}}
		\put(1,5){\makebox[0pt]{\scriptsize #1}}
		\put(7,4){\makebox[0pt]{\scriptsize #2}}
		\put(13,5){\makebox[0pt]{\scriptsize #3}}
	\end{picture}
}
\newcommand{\Dtriangle}[2]{%
	\setlength{\unitlength}{.08cm}
	\rule[-3\unitlength]{0pt}{18\unitlength}
	\begin{picture}(18,7)(0,3)
		\put(4,4){\circle{2}}
		\put(5,4){\line(1,0){8}}
		\put(14,4){\circle{2}}
		\put(4.4472,4.8944){\line(1,2){4.1056}}
		\put(9,14){\circle{2}}
		\put(13.5528,4.8944){\line(-1,2){4.1056}}
		\put(2,3){\makebox[0pt][r]{\scriptsize #1}}
		\put(9,17){\makebox[0pt]{\scriptsize #1}}
		\put(16,3){\makebox[0pt][l]{\scriptsize #1}}
		\put(6,9){\makebox[0pt][r]{\scriptsize #2}}
		\put(12.5,9){\makebox[0pt][l]{\scriptsize #2}}
		\put(9,1){\makebox[0pt]{\scriptsize #2}}
	\end{picture}
}

\makeatother
\begin{document}

\title{On Nichols algebras over $\slfq$ and $\glfq$}

\author{Sebasti\'an Freyre}
\author{Mat\'{\i}as Gra\~na}
\author{Leandro Vendramin}
\thanks{This work was partially supported by CONICET and ANPCyT (Argentina)}
\email{(sfreyre|matiasg|lvendramin)@dm.uba.ar}
\address{Departamento de Matem\'atica - FCEyN \\
Universidad de Buenos Aires \\
Pab I - Ciudad Universitaria \\
(1428) Buenos Aires - Argentina}

\begin{abstract}
	We compute necessary conditions on Yetter-Drinfeld modules over the groups $\slfq$ and
	$\glfq$ to generate finite dimensional Nichols algebras.  This is a first step towards
	a classification of pointed Hopf algebras with a group of group-likes isomorphic to
	one of these groups.
\end{abstract}

\maketitle

\section{Introduction}
Hopf algebras and their variants appear in many contexts; Quantum Groups \cite{MR934283},
Rational and Logarithmic Conformal Field Theories \cite{MR1903984,MR2030633}, Quantum
Field Theories \cite{brouder-2006} come to mind. This paper is part of the long problem of
classifying pointed Hopf algebras. To this end, the best tool known so far is the so
called ``Lifting procedure'' \cite{MR1913436}. The main ingredient of the Lifting
procedure are Nichols algebras on Yetter-Drinfeld modules over group algebras. The
classification of finite dimensional Nichols algebras over group algebras turns out to be
a hard problem: when one restricts the attention to abelian groups, it includes the
classification of semisimple Lie algebras \cite{ARXIV:0605795}; when one considers
non-abelian groups, only a few (genuine) examples are known \cite{zoo} and we possess a
very limited amount of general results (see \cite{MR1994219} for an account).

In last years, two main directions were pursued to classify finite dimensional pointed
Hopf algebras over non-abelian groups. One of them is to generate Nichols algebras by
means of racks and $2$-cocycles \cite{MR1994219}, getting examples that can be present for
many groups. The other one is to concentrate on certain (families of) groups and classify
which finite dimensional Nichols algebras can appear. In the last vein, let us cite the
recently appeared papers: \cite{ARXIV:0511020,ARXIV:0608701,ARXIV:0702559} (where
symmetric groups are considered). The present paper can be thought of as an analog of
these, but for some finite groups of Lie type: $\glfq$ and $\slfq$. As in those
papers, using the tools in \cite{ARXIV:0605795} we are able to rule out most of the
Yetter-Drinfeld modules over these groups, as they generate infinite dimensional Nichols
algebras.

After applying these techniques, the cases we are left with are difficult and one needs
stronger tools to compute them. As an example, we point out that when $q=5$, the conjugacy
classes $\mathcal{C}_i$ ($i=3,\ldots,6$) in Table~\ref{ta:slqi} are isomorphic as racks to
that of the faces of a dodecahedron, while one of the two classes $\mathcal{C}_8$ is
isomorphic to the rack of the faces of an icosahedron (see \cite{zoo}). The Nichols
algebras generated by these racks with a negative $2$-cocycle are not yet known.

One of the main results of this paper is the following theorem concerning Nichols
algebras over the group $\slfq$:

\begin{thm*}
	If the Nichols algebra $\nic(\mathcal{C})$ associated to a conjugacy class
	$\mathcal{C}$ 
of the group $\slfq$ is finite dimensional then $q$ is odd and $\mathcal{C}$ is
one of the classes $\mathcal{C}_i$ of the table \ref{ta:slqi} for $i=2,5,6,7,8$.
\end{thm*}

Necessary conditions over the representations are explicitly given in propositions
\ref{pr:slc2}, \ref{pr:slc7}, \ref{pr:slc8} and \ref{pr:slc5}.

In section \ref{gl} we analyze Nichols algebras over $\glfq$ and as a consequence we
obtain necessary conditions over the representations to get finite-dimensional Nichols algebras
over this group.

\section{Preliminaries and conventions}

\subsection{Yetter Drinfeld modules over $kG$}
In what follows $G$ is a finite group and $k$ is an algebraically closed field of
characteristic $0$.  Recall that a $kG$-comodule $V$ is just a $G$-graded vector space:
$V=\oplus_{g\in G}V_{g}$, where $V_{g}=\{ v\in V\mid\delta(v)=g\otimes v\}$.  A
Yetter-Drinfeld module over $kG$ is a left $kG$-module and a left $kG$-comodule $V$
satisfying the following compatibility condition
\[
	\delta(g\cdot v)=ghg^{-1}\otimes g\cdot v,\quad
	\text{for } v\in V_h
\]
We denote by $\ydg {kG}$ the category of Yetter-Drinfeld modules over $kG$, the morphisms
of which are the maps of modules and comodules.  This is a braided category, with braiding
given by
\begin{equation}\label{eq:prieq}
	c_{M,N}:M\otimes N\to N\otimes M,\quad
	c(m\otimes n)=gn\otimes m,\quad\text{for }m\in M_g.
\end{equation}

Let $g\in G$. We write $\mathcal{Z}_{g}$ for the centralizer
$\mathcal{Z}_g=\{ x\in G\mid xgx^{-1}=g\}$.  We write $\mathcal{C}_g$ for the conjugacy
class $\mathcal{C}_g=\{xgx^{-1}\mid x\in G\}$. Let $\rho:\mathcal{Z}_{g}\to\gl(W)$ be a
representation. We denote by $V(g,\rho)$ the space $\Ind_{\mathcal{Z}_g}^{G}\rho$, endowed
with the comodule structure $\delta(h\otimes w)=ghg^{-1}\otimes(h\otimes w)$. It is
folklore that $V(g,\rho)$ is irreducible iff $\rho$ is, that $\ydg{kG}$ is semisimple iff
$kG$ is, and that in this case the simple modules are given by $V(g,\rho)$, letting $g$
run over representatives of conjugacy classes and $\rho$ on irreducible representations of
$\mathcal{Z}_g$.  In this work we need only to consider representations of degree $1$
since all the centralizers we have to deal with are abelian.

Let $\mathcal{C}=\mathcal{C}_g=\{ g_{i}\mid i\in I\}$ be the conjugacy class of $g$, with
$g_{i}=x_{i}gx_{i}^{-1}$, and let $v_{i}=x_{i}\otimes1\in V(g,\chi)$, where
$\chi=\rho\in\widehat{\mathcal{Z}_g}$ is a character.  Let $T\subseteq I$ be a subset such
that $g_{i}g_{j}=g_{j}g_{i}$ for all $i,j\in T$.  Let $V_{T}\subseteq V(g,\chi)$ be the
subspace generated by $\{ v_{i}\mid i\in T\}$.
Then the braiding restricted to $V_{T}$ is of \emph{diagonal type}, given by
\begin{equation}\label{eq:trnzadiag}
	c(v_{i}\otimes v_{j})=\mbq_{ij}v_{j}\otimes v_{i},\quad
	\text{where}\ \mbq_{ij}=\chi(x_{j}^{-1}g_{i}x_{j})\in k.
\end{equation}
Indeed, by \eqref{eq:prieq}, if $v_j=x_j\otimes w$, we have
\begin{align*}
c(v_i\otimes v_j) &= g_iv_j\otimes v_i = (x_jx_j^{-1}g_ix_j\otimes w) \otimes v_i \\
	&= (x_j\otimes x_j^{-1}g_ix_j w) \otimes v_i = (x_j \otimes \mbq_{ij}w)\otimes v_i,
\end{align*}
since $x_j^{-1}g_ix_j\in\mathcal{Z}_g$.

We write $\mbq=\mbq_T=(\mbq_{ij})$.
If $T\subseteq I$ and $T'\subseteq\mathcal{C}$, $T'=\{g_i\mid i\in T\}$, then 
we abuse notation by calling $\mbq_{T'}=\mbq_T$ and $V_{T'}=V_T$.

\subsection{Notations}
Throughout the paper, $p$ will stand for a prime number and $q$ for a power of $p$.

If $\rho=\chi$ is a character of $\mathcal{Z}_{g}$, we usually write the Nichols algebra
generated by the Yetter-Drinfeld module $V(g,\rho)$ by $\nic(\mathcal{C},\chi)$ or just
$\nic(\mathcal{C})$ when no confusion can arise.

If $n\in\NN$, we write $\ra n\subset k$ for the set of primitive $n$-th roots of unity in
$k$. The order of an element $h$ in a group will be denoted by $|h|$.

We recall from \cite{ARXIV:0605795} the Dynkin diagram notation for braidings of
diagonal type. If $\{v_1,\cdots,v_\ell\}$ is a basis of the braided vector space
$V$ and
$$c(v_i\otimes v_j)=\mbq_{ij}v_j\otimes v_i,$$
then the Dynkin diagram of $c$ has one vertex for each $i=1,\cdots,\ell$, and vertices $i$
and $j$ are joined by an edge iff $\mbq_{ij}\mbq_{ji}\ne1$. Moreover, vertices and edges
have labels as follows: the vertex $i$ has as label the number $\mbq_{ii}$, while the edge
between $i$ and $j$ (if any) has as label the number $\mbq_{ij}\mbq_{ji}$.

\subsection{Needed lemmas}

\begin{lem}\label{lem:tablar2}
	Let $W$ be a two-dimensional space with a braiding of diagonal type. Assume that the
	Dynkin diagram of $W$ is given by \Dchaintwo{$\zeta$}{$\mu$}{$\zeta$} and suppose that
	$\dim\nic(W)<\infty$.  Then the Dynkin diagram is among the following ones: \newline
	\begin{center}
		\begin{tabular}[h]{|c|c|} \hline
			Dynkin diagram & fixed parameter \\ \hline
			\hspace{1cm}\Dospuntos{$\alpha$}{$\alpha$}\hspace{1cm}  &  \hspace{1cm}$\alpha\in k^\times$\hspace{1cm} \\ \hline
			\Dchaintwol{$\alpha$}{$\alpha^{-1}$}{$\alpha$}  &  $\alpha\neq 1$ \\ \hline
			\Dchaintwol{$-1$}{$\alpha$}{$-1$}  &  $\alpha\neq\pm 1$ \\ \hline
			\Dchaintwol{$-\alpha^{-2}$}{$-\alpha^3$}{$-\alpha^{-2}$}  &  $\alpha\in\ra{12}$ \\ \hline
			\Dchaintwol{$-\alpha^{-2}$}{$\alpha$}{$-\alpha^{-2}$}     &  $\alpha\in\ra{12}$ \\ \hline
		\end{tabular}
	\end{center}
\end{lem}
\begin{proof}
	This follows at once by inspection on \cite[Table 1]{ARXIV:0509145}
\end{proof}

The following corollary appears in different forms in
\cite{ARXIV:0511020,ARXIV:0608701,ARXIV:0702559}. We prove it here for the reader's
convenience.
\begin{cor}\label{cor:potinv}\
	\begin{enumerate}
		\item\label{lem:potencias}
			Assume there exists $x\in G$ such that $xgx^{-1}=g^{n}\ne g$. Let
			$T=\{g,g^{n}\}$ and write $\frac 1n=n^{-1}\pmod{|g|}$. Then
			$\mbq=\mdpd{\alpha}{\alpha^{\frac 1n}}{\alpha^{n}}{\alpha}$, $\alpha=\chi(g)$.
			If $\dim\nic(\mathcal{C},\rho)<\infty$ then $\alpha=-1$ or $\alpha\in\ra3$.
			If moreover $g^{n^2}\neq g$, then $\alpha=-1$.

		\item\label{cor:inverso}
			As a particular case, if $g^{n}=g^{-1}$ and
			$\dim\nic(\mathcal{C},\rho)<\infty$ then $g$ has even order
			and $\mbq=\mdpd{-1}{-1}{-1}{-1}$.
	\end{enumerate}
\end{cor}
\begin{proof}
We consider first the case $g^{n^2}=g$. After using Lemma~\ref{lem:tablar2} we obtain that
$\alpha^{n+\frac{1}{n}}=1$ or $\alpha^{n+\frac{1}{n}+1}=1$.  If $|\alpha|$ divides $n^2+1$
then, since $|\alpha|$ divides $n^2-1$, $\alpha=-1$ ($\alpha=1$ would imply
$\dim\nic(\mathcal{C},\rho)=\infty$). If $|\alpha|\mid n^2+n+1$ then $|\alpha|$ divides
$n+2$ and then $\alpha\in\ra3$. Now we consider the case $g^{n^2}\ne g$, i.e.  $g$, $g^n$
and $g^{n^2}$ are different and they belong to $\mathcal{Z}_g\cap\mathcal{C}_g$. Then, if
$T=\{g,g^n,g^{n^2}\}$ we have
\[
\mbq_T=\left(\begin{array}{ccc}
\alpha & \alpha^{\frac{1}{n}} & \alpha^{\frac{1}{n^{2}}}\\
\alpha^{n} & \alpha & \alpha^{\frac{1}{n}}\\
\alpha^{n^{2}} & \alpha^{n} & \alpha\end{array}\right).
\]
By inspection on \cite[Table~2]{ARXIV:0509145}, the only possibilities for $\mbq_T$ to produce
a finite dimensional Nichols algebra are either
\begin{itemize}
	\item $\alpha=-1$, or
	\item $\alpha^{\frac 1n+n}=1$ and $\alpha^{\frac 1{n^2}+n^2}=1$ (but then
		$\alpha=-1$), or
	\item $\alpha^{\frac 1n+n}=1$ and $\alpha^{\frac 1{n^2}+n^2+1}=1$ (but then
		$\alpha=1$), or
	\item $\alpha^{\frac 1{n^2}+n^2}=1$ and $\alpha^{\frac1n+n+1}=1$ (but then
		$\alpha=1$).
\end{itemize}
This completes the proof.
\end{proof}

\begin{lem}\label{lem:tablar3}
	\begin{enumerate}
		\item Let $W$ be a three-dimensional space with a braiding of diagonal type. Assume that the
			Dynkin diagram of $W$ is 
			\begin{center}
				\Dtriangle{$\alpha$}{$\beta$}
			\end{center}
			and suppose that $\dim\nic(W)<\infty$. Then $\alpha=-1$ and
			$\beta\in\ra3$.
		\item Suppose that the Dynkin diagram contains a cycle of length $\ge 4$. Then
			$\dim\nic(W)=\infty$.
	\end{enumerate}
\end{lem}
\begin{proof}
	The first part follows by inspection on \cite[Table 2]{ARXIV:0509145}, while the
	second is \cite[Lemma 20]{ARXIV:0605795}.
\end{proof}

\section{Nichols algebras over $\slfq$}

In this section, $E=\mathbb{F}_{q^2}$ will be the quadratic extension of $\cfq$ and
$\overline{x}$ will be the Galois conjugate of $x\in E$.  Recall that the order of $\slfq$
is $(q-1)q(q+1)$. The conjugacy classes of $\slfq$ ar given in Tables \ref{ta:slqp} and
\ref{ta:slqi} (see \cite[\S 5.2]{MR1153249} for $q$ odd and \cite{MR0369497} for $q$ even).

\subsection{The case $q$ even }
There are $q+1$ conjugacy classes divided in $4$ types:
\begin{longtable}{|c|c|c|c|c|}
\caption{Conjugacy classes in $\slfq$, $q$ even.}\label{ta:slqp} \\ \hline 
Type& Representative& Size& Number& Centralizer\tabularnewline
\hline 
$\mathcal{C}_{1}$
	& $I=c_{1}=\mdpd 1{}{}1$
	& $1$
	& $1$
	& $\slfq$\tabularnewline \hline
$\mathcal{C}_{2}$
	& $c_{2}=\mdpd 11{}1$
	& $q^{2}-1$
	& $1$
	& $\cfq$\tabularnewline \hline 
$\mathcal{C}_{3}$
	& $c_{3}(x)=\mdpd x{}{}{x^{-1}}\; (x\ne1)$
	& $q(q+1)$
	& $\frac{(q-2)}{2}$
	& $\cfq^{\times}$\tabularnewline \hline
$\mathcal{C}_{4}$
	& $c_{4}(x)=\mdpd{}11{x+\overline{x}}\; (x\in E\setminus\cfq)$
	& $(q-1)q$
	& $\frac{q}{2}$
	& cyclic\tabularnewline \hline
\end{longtable}

We first consider the case $q=2$. We have $\slfq\simeq\mathbb{S}_{3}$.
Nichols algebras generated by irreducible Yetter-Drinfeld modules over $\mathbb{S}_{3}$
were studied in \cite{MR1714540}.  They are infinite dimensional except for
$\mathfrak{B}(\tau,\sgn)$, which is $12$-dimensional (here $\tau$ is a transposition and
$\sgn$ is the non-trivial character of its centralizer). 

We now consider the case $q>2$.
\begin{prop}\label{pr:sl2q2n}
Let $q=2^n$ for $n\ge 2$.  For each $i=1,\ldots,4$ and for any representation $\chi$ of
the centralizer of $c_{i}$, we have $\dim\mathfrak{B}(\mathcal{C}_{i},\chi)=\infty$.
\end{prop}
\begin{proof}
We consider each class separately. The class $\mathcal{C}_{1}$ gives the trivial braiding.
For $\mathcal{C}_{2}$, we take, for $a\in\cfq^{\times}$, $x_{a}=\mdpd a{}{}{a^{-1}}$.
Then $g_{a}=x_{a}c_{2}x_{a}^{-1}=\mdpd 1{a^{2}}{}1$.  The centralizer of $c_{2}$ is the
abelian group $\cfq$ embedded in $\slfq$ as $\mdpd 1{\cfq}{}1$.
If $\chi:\cfq\to\mathbb{C}$ is a character of the centralizer, we get
$\mbq_{ab}
=\chi(x_{b}^{-1}g_{a}x_{b})=\chi(g_{ab^{-1}})$,
%	=\chi\mdpd 1{a^{2}b^{-2}}{}1,$
whence $\mbq_{ab}\mbq_{ba} = \chi\mdpd 1{a^{2}b^{-2}+a^{-2}b^{2}}{}1$ (we write
$\mbq_{uv}$ for $\mbq_{g_ug_v}$). Since $q=2^n$, $\chi$ takes values on ${\pm1}$. If
$\mbq_{aa}=1$, then $\dim\mathfrak{B}(\mathcal{C}_{2},\chi)=\infty$, so we may assume
$\mathbf{q}_{aa}=-1$.

If there exists $a\in\cfq\setminus\{0,1\}$ such that 
$\chi\mdpd 1{a^2+a^{-2}}{}1=-1$, then we get
\[
\mbq_{1,a}\mbq_{a,1}=\mbq_{a,a^2}\mbq_{a^2,a}=\cdots=\mbq_{a^m,1}\mbq_{1,a^m}=-1,
\]
where $a^{m+1}=1$. This implies that the space $V_T$ contains a cycle of length $m\ge 3$
with edges labelled by $-1$, whence $\nic(\mathcal{C}_2,\chi)$ is infinite dimensional by
Lemma~\ref{lem:tablar3}.  Assume then that $\chi\mdpd 1{a^2+a^{-2}}{}1=1$ for all $a\neq
0,1$. But then, for all $x$ in the subgroup generated by the elements
of the form $a^2+a^{-2}$ ($a\neq 0,1$), we get $\chi\mdpd 1x{}1=1$. Take now any $a\in\cfq\setminus\{0,1\}$, and let
$r$ be the order of $a^2+a^{-2}$. Since $r$ is odd, $1=(a^2+a^{-2})^r$ is in the subgroup generated by
elements $a^{2i}+a^{-2i}$ for $0<i\le r$, but this contradicts the fact that
$\mbq_{aa}=\chi\mdpd 11{}1=-1$.

For the class $\mathcal{C}_3$, the centralizer of $c_{3}(x)$ is the cyclic group
$\cfq^{\times}$.  It is easy to see that $c_{3}(x)$ and $c_{3}(x^{-1})=c_{3}(x)^{-1}$ are
conjugate, whence they both have the same odd order.  Then by Corollary \ref{cor:potinv}
\eqref{cor:inverso}, $\dim\mathfrak{B}(\mathcal{C}_{3})=\infty$.

Finally, for the class $\mathcal{C}_4$, we have that the centralizer of $c_{4}(x)$ is a
subgroup of the cyclic group $E^{\times}$, hence it is cyclic.  Again, it is easy to see
that both $c_{4}(x)$ and $c_{4}(x)^{-1}=\mdpd {x+\overline{x}}{1}{1}{}$ are conjugate and
they have odd order.  Then, by Corollary \ref{cor:potinv} \eqref{cor:inverso},
$\dim\mathfrak{B}(\mathcal{C}_{4})=\infty$.
\end{proof}

\subsection{The case $q$ odd.}

There are $q+4$ conjugacy classes divided in $8$ types, displayed in Table \ref{ta:slqi}.

We first consider the case $q=3$. The conjugacy class $\mathcal{C}_1$ in Table
\ref{ta:slqi} gives infinite dimensional Nichols algebras. Class $\mathcal{C}_2$ gives
either exterior algebras or infinite-dimensional algebras. Classes $\mathcal{C}_3$,
$\mathcal{C}_4$, $\mathcal{C}_5$ and $\mathcal{C}_6$ give the rack of order $4$ associated
to the vertices of a tetrahedron (see \cite{MR1800709}). Class $\mathcal{C}_7$ is not
present in this case, and class $\mathcal{C}_8$ gives an infinite dimensional Nichols
algebra of diagonal type, after changing the basis in a similar way as in \cite[Remark
5.2.1]{MR1800709}.

\begin{longtable}{|c|c|c|c|c|}
\caption{Conjugacy classes in $\slfq$, $q$ odd.\label{ta:slqi}} \\ \hline 
Type & Representative & Size & Number & Centralizer \\ \hline 
$\mathcal{C}_{1}$ 
& $c_{1}=\left(\begin{array}{cc} 1 \\ & 1\end{array}\right)$
& $1$
& $1$
& $\slfq$ \\
\hline
$\mathcal{C}_{2}$
& $c_{2}=\left(\begin{array}{cc} -1\\ & -1\end{array}\right)$
& $1$
& $1$
& $\slfq$ \\
\hline
$\mathcal{C}_{3}$
& $c_{3}=\left(\begin{array}{cc} 1 & 1\\ & 1\end{array}\right)$
& $\frac{q^{2}-1}{2}$
& $1$
& $\ZZ_2\times\cfq$ \\
\hline 
$\mathcal{C}_{4}$
& $c_{4}=\left(\begin{array}{cc} 1 & x\\ & 1\end{array}\right)$ $(\sqrt x\not\in\cfq)$
& $\frac{q^{2}-1}{2}$
& $1$
& $\ZZ_2\times\cfq$ \\
\hline 
$\mathcal{C}_{5}$
& $c_{5}=\left(\begin{array}{cc} -1 & 1\\ & -1\end{array}\right)$
& $\frac{q^{2}-1}{2}$
& $1$
& $\ZZ_2\times\cfq$ \\
\hline 
$\mathcal{C}_{6}$
& $c_{6}=\left(\begin{array}{cc} -1 & x\\ & -1\end{array}\right)$ $(\sqrt x\not\in\cfq)$
& $\frac{q^{2}-1}{2}$
& $1$
& $\ZZ_2\times\cfq$ \\
\hline 
$\mathcal{C}_{7}$
& $c_{7}(x)=\mdpd x{}{}{x^{-1}}$ $(x\neq \pm 1)$ 
& $q(q+1)$
& $\frac{(q-3)}{2}$
& $\cfq^\times$ \\
\hline
$\mathcal{C}_{8}$
& $c_{8}(x)=\left(\begin{array}{cc} & -1\\ 1 & x+\overline{x}\end{array}\right)$ $(x\in E\setminus\cfq)$
& $(q-1)q$
& $\frac{q-1}{2}$
& cyclic \\
\hline
\end{longtable}

We now consider the case $q>3$.
\begin{rem}\label{rm:sti}
	Conjugation by $\mdpd x{}{}1\in\glfq$ gives an automorphism of $\slfq$ which
	switches $\mathcal{C}_3$ with $\mathcal{C}_4$, and $\mathcal{C}_5$ with
	$\mathcal{C}_6$. Hence, when classifying the irreducible Yetter-Drinfeld modules which
	produce finite dimensional Nichols algebras over $\mathcal{C}_3$, one automatically
	also classifies those over $\mathcal{C}_4$. Ditto for $\mathcal{C}_5$ and
	$\mathcal{C}_6$.

	Furthermore, if one is only interested in the conjugacy classes \emph{as racks}, it
	turns out that $\mathcal{C}_i$ ($i=3,\ldots,6$) are all isomorphic.  Indeed, consider
	the projection $\slfq\to\pslfq$; it is injective when restricted to these classes.
	Further, it maps $\mathcal{C}_3$ and $\mathcal{C}_5$ to the same conjugacy class in
	$\pslfq$ when $-1$ is a square in $\cfq$, and it maps $\mathcal{C}_3$ and
	$\mathcal{C}_6$ to the same conjugacy class when $-1$ is not a square in $\cfq$.
\end{rem}

\begin{prop}\label{pr:slc2}
	$\dim\nic{(\mathcal{C}_1)}=\infty$. If $\dim\nic{(\mathcal{C}_2)}<\infty$, then
	$\nic({C}_2)$ is the exterior algebra.
\end{prop}
\begin{proof}
	The first statement is clear since in this case the braiding is the usual flip
	$x\otimes y\mapsto y\otimes x$, and therefore the Nichols algebra is the symmetric
	algebra.
	In the second case, if $\chi(c_2)=1$ we get again an infinite dimensional Nichols
	algebra. Thus, we must have $\chi(c_2)=-1$, and the braiding is given by
	$x\otimes y\mapsto -y\otimes x$, wich yields an exterior algebra.
\end{proof}

\begin{prop}\label{pr:slc7}
If $\dim\mathfrak{B}(\mathcal{C}_{7})<\infty$, then $\mbq_T=\mdpd{-1}{-1}{-1}{-1}$ and $x$ has
even order, where $T=\{c_7(x),c_7(x^{-1})\}$ (notice that $c_7(x^{-1})$ and $c_7(x)$ are
conjugate).
\end{prop}
\begin{proof}
	It follows from Corollary \ref{cor:potinv} \eqref{cor:inverso}.
\end{proof}
\begin{prop}\label{pr:slc8}
	Let $x\in E\setminus\cfq$. Then $c_8(x)$ and $c_8(x)^{-1}$ are conjugate, and we take
	$T=\{c_8(x),c_8(x)^{-1}\}$. If $\dim\mathfrak{B}(\mathcal{C}_{8})<\infty$, then
	$\mbq_T=\mdpd{-1}{-1}{-1}{-1}$ (in particular, $c_8(x)$ has even order).
\end{prop}
\begin{proof}
	The proof follows from Corollary \ref{cor:potinv} \eqref{cor:inverso}, once we see that
	$c_8(x)$ and $c_8(x)^{-1}$ are conjugate. To this end, we consider
	the quadratic form
	$$\phi(u,v)=u^{2}+uv(x+\overline{x})+v^{2}$$
	over $\cfq$. It is easy to see
	that this form has rank $2$. Then, by \cite[Proposition 4 (1.7)]{MR0344216},
	it represents all elements of $\cfq^{\times}$ and
	in particular it represents $-1$. Thus, there exist $a,c\in\cfq$ such that
	$a^{2}+ac(x+\overline{x})+c^{2}=-1$. This implies that
	$\mdpd a{c+a(x+\overline x)}c{-a}\in\slfq$.  Furthermore,
	\[
		\mdpd a{c+a(x+\overline{x})} c{-a}
			\mdpd {}{-1} 1{x+\overline{x}}
			{\mdpd a{c+a(x+\overline{x})} c{-a}}^{-1}
		=\mdpd {x+\overline{x}}1 {-1}{}
	\]
	finishing the proof.
\end{proof}

\begin{prop}\label{pr:slc3}
	$\dim\mathfrak{B}(\mathcal{C}_{3})=\infty$. 
	The same holds for $\mathfrak{B}(\mathcal{C}_{4})$. 
\end{prop}
\begin{proof}[Proof (case $q\neq 3^{2n+1}$)]
	Let $\alpha=\chi(c_3)$ ($\chi\in\widehat{\mathbb{Z}_2}\times\widehat{\cfq}$).  First of
	all note that $|c_3|=p$. Thus, if $\dim\mathfrak{B}(\mathcal{C}_{3})<\infty$, we may
	suppose that $|\alpha|=p$ (if $\alpha=1$ then $\nic(\mathcal{C}_3)$ is infinite
	dimensional). Let us consider the case $q\equiv1\pmod 4$. In this case there exists an
	element $a\in\mathbb{F}_{q}$ such that $a^{2}=-1$. Then
	\[
	\mdpd a{}{}{a^{-1}}
	\mdpd 11{}1
	\mdpd {a^{-1}}{}{}a
	=\mdpd 1{-1}{}1
	= {\mdpd 11{}1}^{-1}.
	\]
	Now the Proposition follows from Corollary \ref{cor:potinv} \eqref{cor:inverso}, since
	$c_3$ must have an even order.

	If $q\equiv3\pmod 4$ we have (since $p\neq 3$),
	\[
	c_3\neq c_3^{4}=\mdpd 14{}1=\mdpd 2{}{}{2^{-1}}\mdpd 11{}1\mdpd{2^{-1}}{}{}2
	\]
	and then, by Corollary~\ref{cor:potinv} \eqref{lem:potencias}, we get 
	$\dim\mathfrak{B}(\mathcal{C}_{3})=\infty$.

	The statement about $\mathcal{C}_4$ follows from Remark~\ref{rm:sti}.
\end{proof}

We consider now the case $q=3^{2n+1}$. For this, we need Lemma~\ref{lem:lematec} and
Corollary~\ref{co:snl} below, which we could not find in the literature.
Let $\varphi$ denote the Euler function,
$\varphi(n)=\#\{m\in\mathbb{N}\;|\;m<n\text{ and $m$ is coprime to }n\}$.
\begin{lem}\label{lem:lematec}
	If $n\in\mathbb{N}$ is such that $3\nmid n$ and $4\nmid n$,
	then $\varphi(n)>(\frac n2)^{\frac{3}{4}}$.
\end{lem}
\begin{proof}
	If $n$ is even, $n=2p_{1}^{r_{1}}p_{2}^{r_{2}}\cdots p_{N}^{r_{N}}$, where
	$5\le p_{1}<p_{2}<\dots<p_{N}$ are prime numbers. Then
	\[
	\frac{\varphi(n)^{2}}{n}
		=\frac{1}{2}\prod_{i}\frac{(p_{i}-1)^{2}p_{i}^{2r_{i}-2}}{p_{i}^{r_{i}}}
		=\frac{1}{2}\prod_{i}(p_{i}-1)^{2}p_{i}^{r_{i}-2}.
	\]
	Since $(p_{i}-1)^{2}>p_{i}^{3/2}\geq p_{i}^{2-\frac{r_{i}}{2}}$
	for $r_{i}\geq1$, we have
	\[
		\frac{\varphi(n)^{2}}{n}>\frac{1}{2}\prod p_{i}^{\frac{r_{i}}{2}}
		=\frac{\sqrt n}{2\sqrt{2}}.
	\]
	If $n$ is odd,
	$\varphi(n)=\varphi(2n)>n^{\frac{3}{4}}>(\frac{n}{2})^{\frac{3}{4}}$.
\end{proof}

\begin{cor}\label{co:snl}
	$\frac{\varphi(3^{p}-1)}{p}>3^{\frac{p-1}{2}}$ for all odd prime
	number $p$.
\end{cor}
\begin{proof}
	If $p=3$, $5$ or $7$, it follows immediately. Otherwise, by the Lemma,
\[
\frac{\varphi(3^p-1)}{p}>\frac{(3^p-1)^{3/4}}{2^{3/4}p}
	>\frac{3^{3p/4}}{2p} > 3^{\frac{p-1}{2}}.
\]
\end{proof}

We finish now with the proof of Proposition~\ref{pr:slc3}.
\begin{proof}[End of proof of Proposition \ref{pr:slc3} (case $q=3^{2n+1}$)]
	Assume first that $r=2n+1$ is a prime number. For $z\in\cfq$, we denote
	$f_z=f_z(X)\in\mathbb{F}_3[X]$ the minimal polynomial of $z$. We consider the sets 
	\begin{align*}
		\mathcal{I}
			&=\left\{f_{x^2}\in\mathbb{F}_3[X]\;|\;
				x^2\in\mathbb{F}_{3^r}\setminus\mathbb{F}_{3},\,|x^2|=\frac{3^r-1}{2}\right\},
				\text{and} \\
		\mathcal{S}
			&=\left\{X^r+a_{r-1}X^{r-1}+\dots+a_1 X-1\in\mathbb{F}_3[X]\;|\;
			\forall i=1,\ldots,r,\,a_{r-i}+a_i=0\right\}.
	\end{align*}
	Note that $\#\mathcal{I}=\frac{\varphi(\frac{3^r-1}{2})}r=\frac{\varphi(3^r-1)}{r}$
	(because all minimal polynomials of elements in $\mathbb{F}_{3^r}\setminus\mathbb{F}_3$ have
	degree $r$, since $r$ is a prime number).  We have that
	$\#\mathcal{S}=3^\frac{r-1}{2}<\#\mathcal{I}$ by Corollary \ref{co:snl}.
	Then there exists $f\in\mathcal{I}\setminus\mathcal{S}$, and we chose such
	an $f$.  Since $f(0)=(-1)^r\cdot\norm(x^2)=-1$, in order for $f$ not to be
	in $\mathcal{S}$, there exists an $i$ such that $a_{r-i}+a_i\ne0$.

	Consider now $x\in \mathbb{F}_{3^r}^\times$ of order $\frac{3^r-1}{2}$,
	and let $\textbf{X}=\{x,x^3,\dots,x^{3^{r-1}}\}$ 
	be the orbit of $x$ by the action of the Galois group of the extension
	$\mathbb{F}_{3^r}/\mathbb{F}_3$.
	Then $\prod_{y\in\textbf{X}} y = 1$ (because $x$ is a square).
	Also, if $\emptyset\neq\mathbf{S}\subseteq\mathbf{X}$ and
	$\prod_{y\in\mathbf{S}} y \in \mathbb{F}_{3}^\times$,
	then $\mathbf{S} =\mathbf{X}$. Indeed, 
	$\prod_{y\in\mathbf{X}} y = x^{\sum_{i=0}^{r-1}3^i}= x^\frac{3^r-1}{2}= 1$. 
	Now, if $z=\prod_{y\in\mathbf{S}} y$, and $\mathbf{S}\neq\mathbf{X}$,
	then $z=x^s$ for some $s<\frac{3^r-1}{2}$. 
	Note that $x$ is a square and $s<|x|$, which implies that $x^s\ne\pm1$.
	Let
	\begin{align*}
		A &= \{a^2+a^{-2}\mid
		a\in\mathbb{F}_{3^r}^\times\setminus(\ra1\cup\ra2\cup\ra3\cup\ra4\cup\ra6)\}
		= \{a^2+a^{-2}\mid a\in\mathbb{F}_{3^r}\setminus\mathbb{F}_3\}.
	\end{align*}
	Now,
	\begin{equation*}
		a_i = \sum_{\substack{\mathbf{Y}\subset\mathbf{X} \\ \#\mathbf{Y}=r-i}}
			\prod_{y\in\mathbf{Y}}y,
		\qquad
		a_{r-i} = \sum_{\substack{\mathbf{Y}\subset\mathbf{X} \\ \#\mathbf{Y}=i}}
			\prod_{y\in\mathbf{Y}}y
			=\sum_{\#\mathbf{Y}=r-i}\prod_{y\in\mathbf{Y}}\frac 1y.
	\end{equation*}
	Thus, $a_{r-i}+a_i\in (A)$, the subgroup generated by $A$. Then, $1\in (A)$.

	Let now $\chi$ be a character of the centralizer of $c_3$.  If $\chi\mdpd 1t{}1\neq 1$
	for some $t\in A$, then we get a Dynkin diagram with a cycle of length $\ge 4$, which
	is of infinite type. Otherwise, since $1\in(A)$, we get $\mbq_{11}=1$, which also yields
	an infinite dimensional Nichols algebra by Lemma \ref{lem:tablar3}.

	If $r=2n+1$ is not a prime number, one can repeat the same argument with $r'$, a prime factor of
	$r$, taking $\mathbb{F}_{3^{r'}}$ as a subfield of $\mathbb{F}_{3^r}$.
\end{proof}

\begin{prop}\label{pr:slc5}
	If $p\neq 3$ and $\dim\nic(\mathcal{C}_{5})<\infty$, then $\chi=\sgn\times\varepsilon$,
	where $\sgn$ is the non-trivial representation of $\ZZ_2$ and $\varepsilon$ is the
	trivial representation of $\cfq$. The same statement goes for $\nic(\mathcal{C}_6)$.
\end{prop}
\begin{proof}
	As before, for $a\in\cfq^\times$, let
	$$
	x_a=\mdpd a{}{}{a^{-1}},\quad
	g_a=x_ac_5x_a^{-1}=\mdpd {-1}{a^2}{}{-1}.
	$$
	Then $\mbq_{ab}=\chi\mdpd{-1}{a^2b^{-2}}{}{-1}$, and
	$\mbq_{a1}\mbq_{1a}=\chi\mdpd1{-(a^2+a^{-2})}{}1$. For $t\in\cfq$, let $n_t=\mdpd 1t{}1$.
	Let
	$$A=\{-(a^2+a^{-2})\;|\;a\in\cfq^\times\setminus\ra1\cup\ra2\cup\ra3\cup\ra4\cup\ra6\}.$$
	If there exists $t\in A$ with $\chi(n_t)\neq 1$, then $\dim\nic(\mathcal{C}_5)=\infty$ by
	Lemma \ref{lem:tablar3}, as in the proof of Proposition \ref{pr:sl2q2n}.  Thus, we may
	suppose that $\chi (n_t)=1$ for all $t\in A$.  Now, notice that $\# A\ge\frac{q-9}{4}$.
	When $q\not\in\{5,25,7,11,13\}$, we get $\# A>\frac{q}{p}$.  In this case, $A$
	generates (as an abelian group) the whole $\cfq$, which implies the statement.  If
	$q=11$ or $q=13$ an easy computation shows that $A$ generates $\cfq$, whence we are
	done.

	Let $q=25$, $\cfq=\mathbb{F}_5[s]/(s^2-2)$. Let first $a=1+s$, then $|a|=12$ (hence
	$a\not\in\ra1\cup\ra2\cup\ra3\cup\ra4\cup\ra6$), and $-(a^2+a^{-2})=-( (3+2s) + (3+3s))=4$.
	Now, take $a=1+2s$, then $|a|=24$ and $-(a^2+a^{-2})=2s$.  Therefore, $A$ generates
	$\cfq$, and we are done.

	We deal now with the cases $q\in\{5,7\}$.
	Let $q=5$ and $a=2\in\cfq^\times$. With $g_1,g_a$, we get a Dynkin diagram as in
	Lemma~\ref{lem:tablar2}, with $\zeta=\chi\mdpd{-1}1{}{-1}$, and $\mu=\chi\mdpd 12{}1$.
	Since this implies that $|\zeta|$ divides $10$ and $\mu=\zeta^8$, by the Lemma we get
	that $\nic(\mathcal{C}_5)$ is infinite dimensional unless either
	\begin{itemize}
		\item $\zeta\neq 1$ and $\zeta\mu=1$, or
		\item $\zeta\neq 1$ and $\mu=1$, or
		\item $\zeta=-1$.
	\end{itemize}
	The first case is impossible, while in the second and third case we arrive to the
	statement.

	Let $q=7$. Let $a=3$, $\alpha=\chi\mdpd{-1}1{}{-1}$ and $\beta=\chi\mdpd 11{}1$. Then if
	$\beta\neq 1$, we have a Dynkin diagram as in Lemma~\ref{lem:tablar3}, which yields an
	infinite dimensional Nichols algebra (notice that $\beta\not\in\ra3$).  If $\beta=1$,
	as $\alpha^6=\beta$ and $|\alpha|$ divides $14$, this implies $\alpha=-1$, which is
	the statement again.

	Remark~\ref{rm:sti} gives the statement about $\mathcal{C}_6$.
\end{proof}

\begin{rem}
	If $p=3$ one can prove, using the same techniques as in Proposition~\ref{pr:slc3}
	(case $q=3^{2n+1}$), that $\chi=\sgn\times\varepsilon$ when restricted to each
	subfield $\mathbb{F}_{3^r}$ where $r$ is an odd prime.
\end{rem}

\section{Nichols algebras over $\glfq$}\label{gl}
We proceed in this section with the groups $\glfq$. Recall that the order of $\glfq$ is
$(q-1)^2q(q+1)$. Again, $E$ will be the quadratic extension of $\cfq$ and 
$\overline{x}$ will be the Galois conjugate of $x\in E$. Since 
$\mathbf{GL}(2,\mathbb{F}_2)\simeq\mathbb{S}_3$, we consider only the case $q>2$.

There are $q^{2}-1$ conjugacy classes divided in $4$ types:
\begin{longtable}{|c|c|c|c|c|}
\caption{Conjugacy classes in $\glfq$.} \\ \hline 
Type&
Representative&
Size&
Number&
Centralizer\tabularnewline \hline 
$\mathcal{C}_{1}$&
	$c_{1}(x)=\left(\begin{array}{cc}
	x\\
	& x\end{array}\right)$&
	$1$&
	$q-1$&
	$\glfq$\tabularnewline\hline
$\mathcal{C}_{2}$&
	$c_{2}(x)=\left(\begin{array}{cc}
	x & 1\\
	 & x\end{array}\right)$&
	$q^{2}-1$&
	$q-1$&
	$\cfq^\times\times\cfq$\tabularnewline\hline 
$\mathcal{C}_{3}$&
	$c_{3}(x,y)=\left(\begin{array}{cc}
	x & \\
	 & y\end{array}\right)$ $(x\ne y)$ &
	$q(q+1)$&
	$\frac{(q-1)(q-2)}{2}$&
	$\cfq^{\times}\times\cfq^{\times}$\tabularnewline\hline
$\mathcal{C}_{4}$&
	$c_{4}(x)=\left(\begin{array}{cc}
	 & -x\overline{x}\\
	1 & x+\overline{x}\end{array}\right)$ $(x\in E\setminus\cfq)$&
	$(q-1)q$&
	$\frac{q(q-1)}{2}$&
	\text{cyclic} \tabularnewline\hline
\end{longtable}

For the next proposition we need to recall the character table of the non-abelian group
$\glfq$ from \cite[\S 5.2]{MR1153249}.  We know that there are $q^2-1$ irreducible
representations.  The relevant information for our purposes is contained in the following
table:

\begin{longtable}{|c|c|c|c|}
\caption{Representations of $\glfq$ in scalar matrices} \\ \hline 
Representation&
Dimension&
Number&
$c_{1}(x)=
	\left(\begin{array}{cc}
	x\\
	 & x\end{array}\right)$\tabularnewline
	\hline 
$U_{\alpha}$ $(\alpha\in\widehat{\mathbb{F}_{q}^{\times}})$&
	$1$&
	$q-1$&
	$\alpha(x)^{2}$\tabularnewline
	\hline 
$V_{\alpha}$ $(\alpha\in\widehat{\mathbb{F}_{q}^{\times}})$&
	$q$&
	$q-1$&
	$\alpha(x)^{2}$\tabularnewline
	\hline 
$W_{\alpha,\beta}$ $(\alpha,\beta\in\widehat{\mathbb{F}_{q}^{\times}}\mbox{ y }\alpha\ne\beta)$&
	$q+1$&
	$\frac{1}{2}(q-1)(q-2)$&
	$\alpha(x)\beta(x)$\tabularnewline
	\hline 
$X_{\gamma}$ $(\alpha\in\widehat{E_{q}^{\times}})$&
	$q-1$&
	$\frac{1}{2}q(q-1)$&
	$\gamma(x)$\tabularnewline
	\hline
\end{longtable}

\begin{prop}
	If $\dim\nic(\mathcal{C}_1,\chi)<\infty$ then:
	\begin{enumerate}
		\item $-1=\alpha(x)^2$, where $\alpha\in\widehat{\cfq^\times}$; or
		\item $-1=\alpha(x)\beta(x)$, where $\alpha,\beta\in\widehat{\cfq^\times}$; or
		\item $-1=\alpha(x)$, where $\alpha\in\widehat{E^\times}$; or
	\end{enumerate}
\end{prop}
\begin{proof}
	If follows from \cite[Lemma 3.1]{MR1800709}.
\end{proof}

\begin{prop}
	If $\dim\nic(\mathcal{C}_{3})<\infty$ then
	$\mbq_T=\mdpd{\alpha}{\beta}{\beta}{\alpha}$ where $\alpha=\chi(c_{3}(x,y))$,
	$\beta=\chi(c_{3}(y,x))$, and $T=\{c_3(x,y),c_3(y,x)\}$. Furthermore, one of the
	following conditions is satisfied:
	\begin{itemize}
		\item $\beta^{2}=1$ and $\alpha\ne1$
		\item $\beta^{2}\ne1$ and $\alpha\beta^{2}=1$
		\item $\beta^{2}\ne1$ and $\alpha\beta^{2}\ne1$ and $\alpha=-1$
		\item $\beta^{2}\in\mathcal{R}_{12}$ and $\alpha=-\beta^{4}\in\mathcal{R}_{3}$
	\end{itemize}
\end{prop}
\begin{proof}
	Notice that $c_3(x,y)$ and $c_3(y,x)$ are conjugated in $\glfq$ by the involution $\mdpd{}11{}$.
	The result now follows from Lemma~\ref{lem:tablar2}.
\end{proof}

\begin{prop}
	Let $T=\left\{c_4(x),\mdpd{x+\overline{x}}{x\overline{x}}{-1}{}\right\}$.
	If $\dim\nic(\mathcal{C}_{4})<\infty$ then $\mbq_T$ is either $\mdpd{-1}{-1}{-1}{-1}$
	or $\mdpd{\omega}{\omega}{\omega}{\omega}$, where $\omega\in\ra3$.
\end{prop}
\begin{proof}
	Note that the centralizer of $c_{4}(x)$ is the cyclic group
	\[
	\left\{ \mdpd a{-cx\overline{x}}c{a+c(x+\overline{x})}
		:a,c\in\cfq^{\times}\right\},
	\]
	and this group is isomorphic to $E^{\times}$ by
	$a+cx\mapsto\mdpd a{-cx\overline{x}}c{a+c(x+\overline{x})}$.
	Take the involution $\mdpd 1{x+\overline{x}}{}{-1}$ to get
	\[
	\mdpd {x+\overline{x}}{x\overline{x}}{-1}{}
	=\mdpd 1{x+\overline{x}}{}{-1}
	\mdpd {}{-x\overline{x}}{1}{x+\overline{x}}
	\mdpd {1}{x+\overline{x}}{}{-1}
	\in\mathcal{C}_{c_4(x)}.
	\]
	Since $E^{\times}$ is cyclic, there exists $n\in\mathbb{N}$
	such that $c_{4}(x)^{n}=\mdpd{x+\overline{x}}{x\overline{x}}{-1}{}$.
	We get the result by using Corollary~\ref{cor:potinv}\eqref{lem:potencias}.
\end{proof}

\begin{prop}
	If $p=2$ then $\nic(\mathcal{C}_2)$ is infinite dimensional. If $p\neq 2$, $q\neq 9$ and
	$\dim\nic(\mathcal{C}_2)<\infty$ then $\chi\mdpd{x^i}a{}{x^i}=(-1)^i$.
\end{prop}

\begin{proof}
	Similarly to the proof of Proposition~\ref{pr:slc5}, we take $x_a=\mdpd a{}{}1$ and
	we get $g_{a}=\mdpd xa{}x$. Then $\mbq_{ab}=\chi(g_{ab^{-1}})$,
	$\mbq_{a1}\mbq_{1a}=\chi\mdpd{x^2}{x(a+a^{-1})}{}{x^2}$.
	As before, if $\mbq_{1a}\mbq_{a1}\neq 1$ for some
	$a\in\cfq^\times\setminus(\ra1\cup\ra2\cup\ra3)$, then we get a Dynkin diagram of
	infinite type, since
	$$
	\mbq_{1,a}\mbq_{a,1}=\mbq_{a,a^2}\mbq_{a^2,a}=\cdots=\mbq_{a^m,1}\mbq_{1,a^m},\qquad
	|a|=m+1.
	$$
	Assume then that $\forall a\in\cfq^\times\setminus(\ra1\cup\ra2\cup\ra3)$,
	$\mbq_{1a}\mbq_{a1}=1$.  We define $\chi_1\in\widehat{\mathbb{F}_q^\times}$ and
	$\chi_2\in\widehat{\mathbb{F}_q}$ by
	$\chi\mdpd\tau\sigma{}\tau=\chi_1(\tau)\chi_2(\tau^{-1}\sigma)$.  Let
	\[
	A=\left\{ x^{-1}\left(a+a^{-1}\right)\mid 
	a\in\cfq^{\times}\setminus(\ra1\cup\ra2\cup\ra3)\right\}.
	\]
	Thus we assume $\chi_1(x^2)\chi_2(b)=1$ for all $b\in A$. Since the orders of $\cfq^\times$ and
	$\cfq$ are coprime, we get $\chi_2(b)=1$ for all $b\in A$.
	Note that, since $\#(\ra1\cup\ra2\cup\ra3)\leq 4$, $\# A\geq\frac{q-5}{2}$.
	If $\# A>\frac{q}{p}=p^{n-1}$ then $A$ is not contained in any $\cfp$-hyperplane of
	$\cfq$, and thus it generates $\cfq$ as an abelian group.  Therefore, $\chi_2=1$ if
	$p^{n-1}(p-2)>5$. Furthermore, we get $\chi_1(x^2)=1$, from where we conclude with the
	statement.

	We study now the cases $q=2^{n}(n\in\mathbb{N})$ and $q\in\{3,5,7\}$.  Let $q=2^{2n}$.  We
	write $\mathbb{F}_4\subseteq\cfq$,
	$\mathbb{F}_4=\mathbb{F}_2[s]/(s^2+s+1)$. We get $s+s^{-1}=1$, $|s|=3$.
	Let then $\beta:=\mbq_{1s}\mbq_{s1}=\chi\mdpd{x^2}x{}{x^2}$,
	$\alpha:=\chi(c_2)=\chi\mdpd x1{}x$. Notice that $\alpha^{2|x|}=1$ and
	$\beta=\alpha^{|x|+2}$.
	Then either $\beta=1$ (which implies $\alpha=1$ since $|x|$ is odd), or, by using
	Lemma~\ref{lem:tablar3}, $\alpha=-1$, but then $\beta=-1$. In any case, we get an
	infinite-dimensional Nichols algebra.

	The case $q=2^{2n+1}$. Since $2^{2n+1}\equiv 2\pmod 3$, then $\mathcal{R}_{3}=\{1\}$.
	Therefore, $\# A=\frac{2^{n}-2}{2}=2^{n-1}-1$. It is easy to see that $1\not\in A$ and $0\not\in A$.
	Let $(A)$ be the subgroup generated by $A$.  We will prove that $1\in(A)$, which will
	imply $(A)=\cfq$, and we shall be done. Let $r$ be a prime number such that
	$r\mid 2n+1$. If $\xi\in\mathbb{F}_{2^{r}}\setminus\mathbb{F}_{2}$, we denote by
	$f_{\xi}$ its minimal polynomial. On the one hand there exist exactly
	$\frac{2^{r}-2}{r}$ irreducible polynomials in $\mathbb{F}_{2}[X]$ of degree $r$. On
	the other hand,
	\[
		f_{\xi}=X^{r}+a_{r-1}X^{r-1}+\cdots+a_{1}X+1
			=\Pi_{i=1}^{r}(X-\xi_{i}).
	\]
	As done in the proof of Proposition~\ref{pr:slc3} (case $q=3^{2n+1}$), one can
	prove that $a_{r-k}+a_k\in(A)\cap\mathbb{F}_2$.
	It suffices then to prove that there exists $k$ such that $a_{r-k}+a_{k}=1$.
	If $a_{r-k}+a_{k}=0$ for all $k$, then $f_{\xi}\in\mathcal{S}$, where
	$$\mathcal{S}=\{X^{r}+b_{r-1}X^{r-1}+\cdots+b_{1}X+1\mid b_{r-k}=b_{k}\ \forall k\}.$$
	Thence, we have
	$\mathcal{I}=\{f\in\mathbb{F}_{2}[X]\mid f\text{ irreducible}\}\subseteq\mathcal{S}$,
	which is a contradiction because $\#\mathcal{S}=2^{\frac{r-1}{2}}<\frac{2^{r}-2}{r}$
	for any $r>3$.  If $r=3$, $f=(X+1)^3\in\mathcal{S}\setminus\mathcal{I}$, and since
	$\#\mathcal{S}= \#\mathcal{I}$, we also get
	$\mathcal{I}\setminus\mathcal{S}\neq\emptyset$.  Then $\#(A)\geq2^{r-1}+1$ and then
	$(A)=\mathbb{F}_{2^{r}}$.

	The case $q=3$. If $x=1$, we use Corollary~\ref{cor:potinv} \eqref{cor:inverso} to see
	that $\dim\nic(\mathcal{C}_2)=\infty$ for any representation (notice that $|c_2(x)|=3$). If
	$x=-1$, $c_2(x)$ generates its centralizer and again by
	Corollary~\ref{cor:potinv}~\eqref{cor:inverso} we get $\chi(c_2(x))=-1$.

	\smallskip
	The case $q=5$. Note that $|c_2(x)|=5|x|$.  If $x=1$ then
	$\nic(\mathcal{C}_2)$ is infinite dimensional by Corollary~\ref{cor:potinv}~\eqref{cor:inverso}.
	Assume now that $x=2$ or $x=3$. We take
	$a=2$ and then $a+a^{-1}=0$. Let 
	\begin{align*}
		\alpha &= \chi(g)=\mbq_{11}=\mbq_{22}=\mbq_{33}=\mbq_{44}, \\
		\beta  &= \mbq_{12}\mbq_{21}=\mbq_{24}\mbq_{42}
			=\mbq_{34}\mbq_{43}=\mbq_{31}\mbq_{13}=\chi\mdpd4{}{}4=\pm 1\text{ and} \\
		\gamma &= \mbq_{23}\mbq_{32}=\mbq_{14}\mbq_{41}=\chi\mdpd41{}4.
	\end{align*}
	If $\beta=-1$ we get a length-four cycle, thence $\beta=1$. Note that
	$\alpha^{10}=1$ and then $|\alpha|=2,5,10$. But $\alpha^{-2}=\gamma$. If $\gamma=1$
	then $\alpha=-1$ and the representation is completely determined as in the statement.
	If $\gamma\ne1$ we have a Dynkin diagram as in Lemma~\ref{lem:tablar2} with
	$\zeta=\alpha$ and $\mu=\alpha^{-2}$, which is of infinite type.  If $x=4$, take again
	$a=2$. We have that $\mbq_{12}\mbq_{21}=\cdots=\beta=\chi\mdpd1{}{}1=1$ and that
	$\gamma=\mbq_{23}\mbq_{32}=\mbq_{14}\mbq_{41}=\chi\mdpd12{}1$. Again,
	$\alpha^{-2}=\gamma$, and we get the statement as for $x=2,3$.

	The case $q=7$. Note that $|c_2(x)|=7|x|$. The case $x=1$ follows as before from
	Corollary~\ref{cor:potinv}~\eqref{cor:inverso}.  The cases $x=3$ and $x=5$ follow by
	taking $a=3$. In fact, let $\alpha=\mbq_{11}=\chi(c_2(x))$ and
	$\beta=\mbq_{13}\mbq_{31}=\chi\mdpd23{}2$. If $\beta\neq 1$, we get a Dynkin diagram
	with a cycle of length $6$. Then $\beta=1$, but $\beta=\alpha^8$ and since $|\alpha|$
	divides $42$, we must have $\alpha=\pm 1$. If $\alpha=1$ we again have an infinite
	dimensional Nichols algebra, whence $\alpha=-1$. But $c_2(x)$ generates the
	centralizer, and we get the statement. The cases $x=2$ and $x=4$ are easier: we take
	as before $a=3$, $\alpha=\mbq_{11}$, $\beta=\mbq_{13}\mbq_{31}$. Then we get that
	$\alpha^8=\beta=\pm 1$ but $|\alpha|$ divides $21$.  If $x=6$ again we use $a=3$,
	$\alpha=\mbq_{11}=\chi(c_2(x))$, $\beta=\mbq_{13}\mbq_{31}=\chi\mdpd1{-1}{}1$. As
	before, we must have $\alpha^{36}=\beta=1$ but $|\alpha|$ divides $14$, which implies
	the result.
\end{proof}

\begin{acknowledgement*}
We benefited from discussions with N.~Andruskiewitsch, F.~Fan\-ti\-no, C.~S\'anchez and
A.~Pacetti, to whom we thank.
We used GAP \cite{GAP} to do some of the computations.
\end{acknowledgement*}

%\bibliographystyle{alphanum}
%\bibliography{referencias.bib}

\begin{thebibliography}{Gra2}

\bibitem[AF1]{ARXIV:0702559}
Nicol{\'a}s Andruskiewitsch and Fernando Fantino.
\newblock On pointed {H}opf algebras associated with alternating and dihedral
  groups.
\newblock math.QA/0702559, 2007,
\newblock To appear in \emph{Rev. Uni{\'o}n Mat. Argent.} 48(1).

\bibitem[AF2]{ARXIV:0608701}
Nicol{\'a}s Andruskiewitsch and Fernando Fantino.
\newblock On pointed {H}opf algebras associated with unmixed conjugacy classes
	in $\mathbb{S}_m$.
\newblock {\em J. Math. Phys.}, 48(3):033502, 26, 2007.

\bibitem[AG1]{MR1714540}
Nicol{\'a}s Andruskiewitsch and Mat{\'{\i}}as Gra{\~n}a.
\newblock Braided {H}opf algebras over non-abelian finite groups.
\newblock {\em Bol. Acad. Nac. Cienc. (C\'ordoba)}, 63:45--78, 1999.
\newblock Colloquium on Operator Algebras and Quantum Groups (Spanish)
  (Vaquer\'\i as, 1997).

\bibitem[AG2]{MR1994219}
Nicol{\'a}s Andruskiewitsch and Mat{\'{\i}}as Gra{\~n}a.
\newblock From racks to pointed {H}opf algebras.
\newblock {\em Adv. Math.}, 178(2):177--243, 2003.

\bibitem[AS]{MR1913436}
Nicol{\'a}s Andruskiewitsch and Hans-J{\"u}rgen Schneider.
\newblock Pointed {H}opf algebras.
\newblock In {\em New directions in {H}opf algebras}, volume~43 of {\em Math.
  Sci. Res. Inst. Publ.}, pages 1--68. Cambridge Univ. Press, Cambridge, 2002.

\bibitem[AZ]{ARXIV:0511020}
Nicol{\'a}s Andruskiewitsch and Shouchuan Zhang.
\newblock On pointed {H}opf algebras associated to some conjugacy classes in
	$\mathbb{S}_n$.
\newblock {\em Proc. Amer. Math. Soc.}, 135(9):2723--2731 (electronic), 2007.

\bibitem[Bro]{brouder-2006}
Christian Brouder.
\newblock Quantum field theory meets {H}opf algebra, 2006.
\newblock hep-th/0611153.

\bibitem[Dri]{MR934283}
V.~G. Drinfel{\cprime}d.
\newblock Quantum groups.
\newblock In {\em Proceedings of the International Congress of Mathematicians,
  Vol. 1, 2 (Berkeley, Calif., 1986)}, pages 798--820, Providence, RI, 1987.
  Amer. Math. Soc.

\bibitem[FH]{MR1153249}
William Fulton and Joe Harris.
\newblock {\em Representation theory}, volume 129 of {\em Graduate Texts in
  Mathematics}.
\newblock Springer-Verlag, New York, 1991.
\newblock A first course, Readings in Mathematics.

\bibitem[Gab]{MR2030633}
Matthias~R. Gaberdiel.
\newblock An algebraic approach to logarithmic conformal field theory.
\newblock {\em Internat. J. Modern Phys. A}, 18(25):4593--4638, 2003.

\bibitem[GAP]{GAP}
The~GAP {G}roup.
\newblock Gap -- Groups, Algorithms, and Programming, version 4.4.9,
  \texttt{http://www.gap-system.org}, 2006.

\bibitem[Gra1]{zoo}
Mat\'{\i}as Gra{\~n}a.
\newblock Nichols algebras of nonabelian group type.
\newblock Web page at \texttt{http://mate.dm.uba.ar/$\sim$matiasg/zoo.html}.

\bibitem[Gra2]{MR1800709}
Mat{\'{\i}}as Gra{\~n}a.
\newblock On {N}ichols algebras of low dimension.
\newblock In {\em New trends in {H}opf algebra theory (La Falda, 1999)}, volume
  267 of {\em Contemp. Math.}, pages 111--134. Amer. Math. Soc., Providence,
  RI, 2000.

\bibitem[Hec1]{ARXIV:0509145}
Istv{\'a}n Heckenberger.
\newblock Classification of arithmetic root systems of rank 3.
\newblock {\em Actas of XVI Colloquium Latinoamericano de \'Algebra}, 2005.

\bibitem[Hec2]{ARXIV:0605795}
Istv{\'a}n Heckenberger.
\newblock Classification of arithmetic root systems, 2006.

\bibitem[PZ]{MR1903984}
Valentina Petkova and Jean-Bernard Zuber.
\newblock Conformal field theories, graphs and quantum algebras.
\newblock In {\em MathPhys odyssey, 2001}, volume~23 of {\em Prog. Math.
  Phys.}, pages 415--435. Birkh\"auser Boston, Boston, MA, 2002.

\bibitem[Ser]{MR0344216}
Jean~Pierre Serre.
\newblock {\em A course in arithmetic}.
\newblock Springer-Verlag, New York, 1973.
\newblock Translated from the French, Graduate Texts in Mathematics, No. 7.

\bibitem[ZN]{MR0369497}
A.~V. Zelevinski{\u\i} and G.~S. Narkunskaja.
\newblock Representations of the group {${\rm SL}(2,\,F\sb{q})$}, where
  {$q=2\sp{n}$}.
\newblock {\em Funkcional. Anal. i Prilo\v zen.}, 8(3):75--76, 1974.

\end{thebibliography}

\def\cprime{$'$}

\end{document}